\documentclass[12pt]{amsart}
\usepackage{amssymb,amsthm}
\usepackage[all]{xy}

\def\silentfootnote#1{{\let\thefootnote\relax\footnotetext{#1}}}

\numberwithin{equation}{section}
\theoremstyle{plain}
\newtheorem{thm}[equation]{Theorem}
\newtheorem{lem}[equation]{Lemma}
\newtheorem{prop}[equation]{Proposition}
\newtheorem{cor}[equation]{Corollary}

\theoremstyle{definition}

\newtheorem{ex}[equation]{Example}

\DeclareMathOperator{\HH}{HH}
\DeclareMathOperator{\coh}{H}
\DeclareMathOperator{\Hom}{Hom}
\DeclareMathOperator{\Ext}{Ext}

\DeclareMathOperator{\End}{End}
\DeclareMathOperator{\Ker}{Ker}
\DeclareMathOperator{\Ima}{Im}
\DeclareMathOperator{\res}{res}

\DeclareMathOperator{\cores}{cor}
\DeclareMathOperator{\id}{id}

\newcommand{\C}{\mathbb{C}}
\newcommand{\kf}{\Bbbk}
\newcommand{\Z}{\mathbb{Z}}

\newcommand{\ba}{\begin{array}}
\newcommand{\ea}{\end{array}}

\title[Hochschild cohomology and Grothendieck rings]
{Products in Hochschild cohomology and Grothendieck rings
of group crossed products}
\author{S.\ J.\ Witherspoon}
\address{Department of Mathematics and Computer Science,
Amherst College, Amherst, Massachusetts 01002}
\email{sjw@cs.amherst.edu}
\date{October 23, 2002}

\begin{document}

\begin{abstract}
We give a general construction of rings
graded by the conjugacy classes of a finite group.
Some examples of our construction are the Hochschild cohomology ring
of a finite group algebra, the Grothendieck ring of the Drinfel'd double
of a group, and the orbifold cohomology ring for a global quotient.
We generalize the first two examples by deriving product formulas
for the Hochschild cohomology ring of a group crossed product and for
the Grothendieck ring of an abelian extension of Hopf algebras.
Our results account for similarities in the
product structures among these examples.

\end{abstract}

\maketitle

\silentfootnote{Research supported by National Security Agency
Grant \#MDA904-01-1-0067.}

\section{Introduction}
Cibils and Solotar first noticed that the {\em cup} product in the Hochschild
cohomology of a finite group algebra is similar to the 
{\em tensor} product
in the category of modules for the Drinfel'd (or quantum) double
of the group (or equivalently of Hopf bimodules for the group algebra)
\cite{cibils97,cibils-solotar97}.
This observation led them to 
conjecture a particular formula for the product in the Hochschild cohomology
ring, which was proven by Siegel and this author \cite{siegel-witherspoon99}.
However, a question remained: 
{\em Why} are these product structures so similar?
Recent work of Bouc provides an answer to this question; they both give
examples of rings constructed in a prescribed way from Green
functors (the cohomology ring of the group, and the Grothendieck ring 
of modules for the group algebra, respectively) \cite{bouc02}.

However there are more examples of cohomology and Grothendieck
rings having similar product structures that do not arise from Bouc's
construction.
These examples have components that are indexed by conjugacy classes of
the group (or more generally by orbits under another group action), 
but that do not come from an underlying Green functor as the components
are not themselves algebras.
The product satisfies a particular formula with respect to this
additive structure, analogous to product formulas for Hochschild
cohomology of the group algebra and for modules of the Drinfel'd double
of the group.
In this paper
we follow the spirit of Bouc's work, recording (in Section 2)
precisely the conditions necessary to recover these
examples.  We generalize Bouc's ring construction in Theorem \ref{invariants},
where we give two potential product formulas, (i) and (ii):
The second (ii) is a direct generalization of Bouc's formula.
The first (i) differs from (ii) by integer scalar factors,
and applies in case the conjugacy class-graded ring is a subring
of invariants of a larger ring (which may not always be the case).
Thus formula (ii) applies to more examples and retains more information
in positive characteristic.

In the remainder of Section 2, we describe several examples of rings
graded by conjugacy classes in which the product structure is known
and coincides with our construction.
In Sections 3 and 4, which are largely self-contained, we
prove product formulas for two further classes of examples.
These examples are all of
independent interest, and justify our construction.
While most of them do not arise from Bouc's construction for Green
functors, they retain many properties of Green functors leading
to their distinctive product structure.

Example \ref{orbifold} is orbifold cohomology.
Given an orbifold with (almost) complex structure, Chen and Ruan
defined an orbifold cohomology ring \cite{chen-ruan00}.
In case the orbifold is a global quotient, that is the quotient of
a manifold by an action of a finite group, Fantechi
and G\"{o}ttsche \cite{fantechi-gottsche01} and Uribe \cite{uribe01}
gave an equivalent definition of its orbifold cohomology ring.
As a consequence of their results, the cup product in the orbifold
cohomology of a global quotient
coincides with that given by our Theorem \ref{invariants}(i). 
If orbifold cohomology were defined more generally for a quotient
of a topological space by a finite group action, 
we expect it would also satisfy the properties
given in Section 2, and hence our Theorem \ref{invariants} would endow
it with a ring structure.
The difference between the two possible products defined in Theorem
\ref{invariants} would be meaningful, as the coefficients may be a ring
of positive characteristic.
We offer the classifying space of a finite group as an example to
illustrate this difference.  Its orbifold cohomology should be the
Hochschild cohomology of the group algebra.

In Section 3, we prove a product formula for
the Hochschild cohomology ring of the crossed product $S\# G$
of an algebra $S$ and a finite group $G$, thus generalizing the product
formula for Hochschild cohomology of a group algebra
\cite[Thm.\ 5.1]{siegel-witherspoon99}.
As a consequence, such a cohomology ring satisfies the properties given
in Section 2, and its product coincides with the product given in our
Theorem \ref{invariants}(ii).
Such crossed products and their (co)homology are featured in many papers:
Lorenz and Cornick gave an additive decomposition of the Hochschild
homology of group-graded rings (of which crossed products are a special case)
into components indexed
by conjugacy classes of the group \cite{cornick95,lorenz92}, generalizing
the well-known case of a group algebra.
Several authors studied the (co)homology of specific types of crossed 
products.
Alev, Farinati, Lambre and Solotar \cite{alev-farinati-lambre-solotar00}
and Alvarez \cite{alvarez02} studied the cohomology of a crossed product
of a Weyl algebra with a finite group, and their calculations were used
by Etingof and Ginzburg \cite{etingof-ginzburg02} to understand
deformations of these crossed products.
One class of such examples arises from an orbifold that is a global quotient.
In this case, $G$ acts on a corresponding commutative algebra $S$ 
of functions, and the (noncommutative)
crossed product ring $S\# G$ (as well as the commutative
ring of invariants $S^G$) is of interest in connection with the orbifold.
See for example Connes' book \cite{connes94}.  
When such a crossed product has a nontrivial twisting
cocycle (see the definitions in Section 3),
the orbifold is said to have discrete torsion, and this case has
been of particular interest.
Caldararu, Giaquinto and this author \cite{caldararu-giaquinto-witherspoon02}
used Hochschild cohomology to find deformations of such a crossed
product arising from an orbifold with discrete torsion that is related
to an example of Vafa and Witten \cite{vafa-witten95}.

A comparison of Example \ref{orbifold} (orbifold cohomology) and
Theorem \ref{productformula} (Hochschild cohomology of a crossed product)
shows that the product structures of these cohomology rings
associated to a given orbifold are very similar.
Our construction in Section 2 of rings graded by conjugacy classes,
of which each of these rings is an example, accounts for this similarity.

In Section 4, we consider representations of 
a finite abelian extension of a Hopf
algebra, that is a Hopf algebra $H$ that is a crossed product
arising from an action of a finite group $L$ on another group $G$.
Such a Hopf algebra is perhaps the simplest possible extension, and
is a fundamental building block for other finite
dimensional Hopf algebras.
(For example all semisimple Hopf algebras
of dimension 16 (over $\C$) have this form \cite{kashina00}.)
In case $H$ is semisimple, all simple $H$-modules were determined by
Kashina, Mason, and Montgomery \cite{kashina-mason-montgomery02};
they are partitioned into classes indexed by the $L$-orbits of $G$.
Here we first
use Clifford theory to extend their result to the nonsemisimple case.
We then prove a formula for the tensor product of any two such
$H$-modules, and show that it coincides with the product given by
our Theorem \ref{invariants}(ii).
This generalizes a formula for the tensor product
of modules for the Drinfel'd double of $G$
(or equivalently of Hopf bimodules of the group algebra \cite{cibils97}).
Our methods work equally well when $H$
is only a {\em quasi} Hopf algebra, and thus apply to the
{\em twisted} Drinfel'd double of $G$ as well.
Our Theorem \ref{invariants} gives 
one explanation of why this formula for the tensor
product of $H$-modules is similar to the product formula for the
Hochschild cohomology ring $\HH^*(\kf L,\kf G)$ of Example \ref{hh*LG}.

The author thanks Alejandro Adem, Don Passman, Arun Ram, and 
Frank Sottile for helpful discussions.

\section{Hochschild constructions}

Let $G$ and $L$ be finite groups with a left action of $L$ on $G$
by automorphisms.
(Often we will let $L=G$ act on $G$ by conjugation.)
Write the action of $x\in L$ on $g\in G$ as ${}^xg$,
and let $L_g$ denote the stabilizer of $g$ in $L$,
that is $$L_g:=\{x\in L\mid {}^xg=g\}.$$

Let $\{A(g)\mid g\in G\}$ be a set of free modules over a commutative
ring $\kf$.
Assume there are $\kf$-linear and $\kf$-bilinear maps
\begin{equation}\label{maps}
c_{g,x}:A(g)\rightarrow A({}^xg) \ \mbox{ and } \
m_{g,h}:A(g)\times A(h)\rightarrow A(gh)
\end{equation}
and that $c_{g,x}$ is an isomorphism of $\kf$-modules, for all $g,h\in G,
x\in L$. Let 
$$A:=\bigoplus_{g\in G}A(g)$$ 
and let $c_x:A\rightarrow A$ be the linear
function which is $c_{g,x}$ on the component $A(g)$.
If $\alpha\in A$, write $\alpha=\sum_{g\in G}\alpha_g$ with $\alpha_g\in A(g)$.
We will also denote by $\alpha_g$ any element of $A(g)\subset A$.
Assume the following properties hold:
\begin{itemize}
\item[(H1)] (group action) $\ c_1=\id$ and $c_x\circ c_y=c_{xy}$ for all
$x,y\in L$, 
\item[(H2)] (compatibility) $c_x\circ m_{g,h}=m_{{}^xg,{}^xh}\circ
(c_x\times c_x)$ for all $g,h\in G, x\in L$, and
\item[(H3)] (unity) there is an element $1\in A(1)$ with $c_x(1)=1$
for all $x\in L$ and $m_{1,g}(1,\alpha_g)=\alpha_g=m_{g,1}(\alpha_g,1)$
for all $g\in G, \alpha_g\in A(g)$.
\end{itemize}
Notice that property (H1) implies there is 
an action of $L$ on the $\kf$-module $A$.
Further, $A$ will be an associative algebra
with an action of $L$ as automorphisms if the following additional
property is satisfied:
\begin{itemize}
\item[(H4)] (associativity) 
$$m_{de,f}(m_{d,e}(\alpha_d,\beta_e),\gamma_f)= 
m_{d,ef}(\alpha_d,m_{e,f}(\beta_e,\gamma_f))$$
for all $d,e,f\in G,
\alpha_d\in A(d),\beta_e\in A(e),\gamma_f\in A(f)$.
\end{itemize}
In this case the product on $A$ is defined componentwise by the maps
$m_{g,h}$, and 
$A$ is a group-graded algebra, graded by the group $G$.
The action of $L$ by automorphisms is compatible with the grading.
Let 
$$A^L:=\{\alpha\in A\mid c_x(\alpha)=\alpha\mbox{ for all }x\in L\},$$
the $\kf$-submodule of $L$-invariants of $A$, also an associative 
algebra in case properties (H1)--(H4) are satisfied.

We would like to weaken property (H4) and still obtain an associative algebra
structure on the invariants $A^L$.
The example of Hochschild cohomology of the group algebra $\kf L$,
with product formula given in \cite[Thm.\ 5.1]{siegel-witherspoon99} and
generalized in \cite[Thm.\ 6.1]{bouc02}, suggests how to do this.
Consider the following property instead:
\begin{itemize}
\item[(H4$'$)] (associativity) 
$$\sum_{(d,e,f)\in T_g}m_{de,f}(m_{d,e}(\alpha_d,\beta_e),\gamma_f)=
\sum_{(d,e,f)\in T_g} m_{d,ef}(\alpha_d,m_{e,f}(\beta_e,\gamma_f)),$$ 
for all $g\in G$, and $\alpha,\beta,\gamma\in A^L$,
where $T_g$ is the set described below.
\end{itemize}
The set $T_g$ is any set of representatives of equivalence classes in
$$\{(d,e,f)\in G\times G\times G\mid def=g\}$$ given by orbits 
under the following group actions.
The stabilizer $L_g$ of $g$ acts on pairs $(de,f)$, that is on the
product of the first two factors, and the last factor.
For each $f$, the stabilizer $L_{gf^{-1}}$ ($=L_{de}$)
acts on pairs $(d,e)$, that is on the first two factors.
We may equally well consider $T_g$ to be a set of
representatives of equivalence classes given by orbits under the action of 
$L_g$ on pairs
$(d,ef)$ and for each $d$, the action of $L_{d^{-1}g}$ ($=L_{ef}$)
on pairs $(e,f)$: The bijection $T_g\rightarrow T_g$ given by
$(d,e,f)\mapsto (gfg^{-1},d,e)$ yields a correspondence of
orbit representatives.
By property (H2),
the choice of representatives does not matter as we are working in
the ring of invariants $A^L$.

Clearly property (H4) implies (H4$'$), but not conversely.
We have already observed part (i) of the following theorem,
and part (ii) describes a product on $A^L$ in case the above weaker
property (H4$'$) holds.

\begin{thm}\label{invariants}
Let $\{A(g)\mid g\in G\}$ be a set of free $k$-modules, with maps
$c_{g,x}$ and $m_{g,h}$ as in (\ref{maps}).
\begin{itemize}
\item[{\rm (i)}] 
If $A$ satisfies (H1)--(H4), then $A$ is an associative ring with
multiplication $m$ given componentwise by $m_{g,h}$ for all $g,h\in G$,
and $A^L$ is a subring under the inherited multiplication.

\item[{\rm (ii)}] 
If $A$ satisfies (H1)--(H3) and (H4$'$), then $A^L$ is an
associative ring with product defined componentwise by
$$(\alpha\cdot \beta)_g=
  \sum_{\substack{(h,k)\in L_g\backslash G\times G\\hk=g}}
      m_{h,k}(\alpha_h,\beta_k)$$
for all $\alpha,\beta\in A^L$ and $g\in G$.\
\end{itemize}
\end{thm}

We remark that if properties
(H1)--(H4) are satisfied, then the $L$-invariants $A^L$ have
two {\em distinct} products given by Theorem \ref{invariants}(i) and (ii).
The product on $A^L$ inherited from $A$ as in (i) is given by
\begin{equation}\label{comparison}
(\alpha\cdot\beta)_g=
  \sum_{\substack{(h,k)\in G\times G\\hk=g}}
m_{h,k}(\alpha_h,\beta_k)=\sum_{\substack{(h,k)\in L_g\backslash G\times G\\
hk=g}}|L_g:L_h\cap L_k| m_{h,k}(\alpha_h,\beta_k),
\end{equation}
where $|L_g:L_h\cap L_k|$ is the index of the subgroup $L_h\cap L_k$ in $L_g$.
This differs from the product (ii) by the integer scalar factors
$|L_g:L_h\cap L_k|$.
Thus the distinction is meaningful in positive characteristic,
where the product (ii) potentially contains more information.

\begin{proof}[Proof of Theorem \ref{invariants}]
The proof of (i) is straightforward.
We will prove (ii).
First note that by property (H2), if $\alpha,\beta\in A^L$, 
then $\alpha\cdot\beta\in A^L$
as well, as the product formula in (ii) yields $(\alpha\cdot\beta)_{{}^xg}
={}^x((\alpha\cdot\beta)_g)$.
By property (H3), $A^L$ has a multiplicative identity.
We will show that (H4$'$) is equivalent to the associativity
of the product defined in (ii).
Let $\alpha,\beta,\gamma\in A^L$ and $g\in G$. Then
$$((\alpha\cdot\beta)\cdot\gamma)_g=
\sum_{\substack{(h,k)\in L_g\backslash (G\times G), \ hk=g\\
(d,e)\in L_h\backslash (G\times G), \ de=h}}m_{h,k}(m_{d,e}(\alpha_d,
\beta_e),\gamma_k),$$
which is the same as the left side of property (H4$'$), whereas
$$(\alpha\cdot(\beta\cdot\gamma))_g=\sum_{\substack
{(h,k)\in L_g\backslash (G\times G), \ hk=g\\(e,f)\in L_k\backslash
(G\times G), \ ef=k}}m_{h,k}(\alpha_h,m_{e,f}(\beta_e,\gamma_f)),$$
the right side of property (H4$'$).
\end{proof}

Let $g_1,\ldots,g_t$ be a set of representatives of orbits of $L$ on $G$.
There is another way to write the product of Theorem \ref{invariants}(ii),
in terms of the additive decomposition
\begin{equation}\label{decomposition}
  A^L=\biggl(\bigoplus_{g\in G}A(g)\biggr)^L \cong
  \bigoplus_{i=1}^t A(g_i)^{L_{g_i}}.
\end{equation}
Write $L_i:=L_{g_i}$ and $\alpha_i:=\alpha_{g_i}$ for $\alpha\in A$.

\begin{cor}\label{product}
Let $\alpha_i\in A(g_i)^{L_{i}}$ and $\beta_j\in A(g_j)^{L_j}$ 
in the additive isomorphism
(\ref{decomposition}). Then under the product of Theorem \ref{invariants}(ii),
$$\alpha_i\cdot \beta_j = \sum_{x\in D}m_{{}^yg_i,{}^{yx}g_j}
  ({}^y\alpha_i,{}^{yx}\beta_j),$$
where $D$ is a set of representatives of double cosets $L_i\backslash
L/L_j$, and $k=k(x)$ and $y=y(x)$ are chosen so that 
${}^yg_i{}^{yx}g_j=g_k$.
\end{cor}

\begin{proof}
The elements $\alpha_i$ and $\beta_j$ correspond in $A^L$ to
$\sum_{y\in L/L_i} {}^y\alpha_i$ and $\sum_{z\in L/L_j}{}^z\beta_j$,
respectively.
By Theorem \ref{invariants}(ii),
$$\biggl(\Bigl(\sum_{y\in L/L_i}{}^y\alpha_i\Bigr)\cdot\Bigl(
\sum_{z\in L/L_j}{}^z\beta_j\Bigr)\biggr)_{g_k}=
\sum_{\substack{({}^yg_i,{}^zg_j)\in L_k\backslash (G\times G)\\
{}^yg_i{}^zg_j=g_k}}m_{{}^yg_i,{}^zg_j}({}^y\alpha_i,{}^z\beta_j).$$
Letting $x=y^{-1}z$ so that $z=yx$, this looks like the sum in the
corollary, other than the set over which the sum is taken.
Noting that the choice of $y=y(x)$ in the corollary is unique only
up to multiplication by an element of $L_k$, and similarly in the
above sum, $y,z$ are unique only up to multiplication by elements
of $L_i,L_j$, respectively, we see that the two sums are the same.
\end{proof}

We give several examples in the remainder of this section.
The first example is just the group algebra.

\begin{ex}\label{groupalgebra} (Group algebra.)
Let $A(g):=\kf g$ for each $g\in G$, that is $A(g)$ is the one-dimensional
vector space with basis $g$.
Define linear maps $m_{g,h}:\kf g\times \kf h\rightarrow \kf gh$ by 
$m_{g,h}(g,h)=gh$. Then $A$ satisfies properties (H1)--(H4), and 
$A=\kf G$ under the product of Theorem \ref{invariants}(i).
The ring of invariants $A^L$ is $(\kf G)^L$.
In particular, if $L=G$, then $A^L$ is the center of $\kf G$.
\end{ex}

The invariant subring $(\kf G)^L$ of the above example 
is isomorphic to $\HH^0(\kf L, \kf G)$,
the degree 0 component of the Hochschild cohomology ring
of $\kf L$ with coefficients in $\kf G$.
The next example is thus a generalization.
The definition of Hochschild cohomology is given in Section 3;
see also \cite{benson91,evens91} for group cohomology.

\begin{ex}\label{hh*LG} (Hochschild cohomology.)
For each $g\in G$, let $A(g):=H^*(L_g)=\Ext^*_{\kf L_g}(\kf,\kf)$ where 
$L_g$ acts trivially on the coefficient ring $\kf$.
The conjugation maps are given by conjugation of group cohomology rings, and 
\begin{equation}\label{mgh}
m_{g,h}(\alpha_g,\beta_h):=\cores^{L_{gh}}_{L_g\cap L_h}(\res^{L_g}
_{L_g\cap L_h}\alpha_g\smile \res^{L_h}_{L_g\cap L_h}\beta_h),
\end{equation}
where $\res$ and $\cores$ denote the standard restriction and corestriction
maps of group cohomology.
Then (H1)--(H3) are well-known properties of group cohomology.
A comparison of \cite[Lem.\ 4.2 and Thm.\ 5.1]{siegel-witherspoon99}
with Corollary \ref{product} above shows that (H4$'$) holds, and that
under the product of Theorem \ref{invariants}(ii), 
$A^L\cong\HH^*(\kf L,\kf G)$ as rings.
\end{ex}

In the above example, the difference between the products (i) and (ii)
of Theorem \ref{invariants} is important:
If there were a definition of $m_{g,h}$ satisfying
(H4) instead of (H4$'$), the subring of invariants $A^L$ of $A$ would in 
general {\em not} be isomorphic to $\HH^*(\kf L,\kf G)$.
For example, if $\kf$ is a field and $n>0$, then $\HH^n(\kf L,\kf G)\neq 0$ 
only if the characteristic of $\kf$ divides $|L|$.
In this case, the integer scalar factor $|L_{gh}:L_g\cap L_h|$ of
(\ref{comparison}) may often be 0 even when the corresponding product
in Hochschild cohomology is nonzero.

The next example, due to Bouc, is a generalization of Example
\ref{hh*LG} arising from any Green functor.
This was the motivating example for our present work.

\begin{ex} \label{greenfunctor}
(Green functor.)
Let $L=G$ and let $A(-)$ be any Green functor for $G$ over $\kf$.
That is, to each subgroup $H<G$, there is assigned a $\kf$-algebra $A(H)$,
together with conjugation maps $c_x:A(H)\rightarrow A({}^xH)$ 
(for all $x\in G$ and $H<G$),
restriction maps $r^H_K:A(H)\rightarrow A(K)$ and transfer maps
$t^H_K:A(K)\rightarrow A(H)$ (for all $K<H<G$). 
These maps are required to satisfy certain properties
(see for example \cite{thevenaz??}).
Let $A(g):=A(C(g))$ where $C(g)=G_g$ is the centralizer of $g$ in $G$, and
$$m_{g,h}(\alpha_g,\beta_h):=t^{C(gh)}_{C(g)\cap C(h)}\left
(r^{C(g)}_{C(g)\cap C(h)}
\alpha_g \cdot r^{C(h)}_{C(g)\cap C(h)}\beta_h\right).$$
Then (H1)--(H3) are standard properties of Green functors.
Property (H4$'$) is equivalent to the associativity of the product on $A^G$
proven by Bouc in \cite[Thm.\ 6.1]{bouc02}.
The ring $A^G$ of our Theorem \ref{invariants}(ii) is precisely 
Bouc's ring $A_G(G)$.
Examples include, among others, the crossed Burnside ring (obtained from the
Burnside ring functor), the Hochschild cohomology ring of $\kf G$
(obtained from the group cohomology functor), and the Grothendieck ring
of the Drinfel'd double of $G$ (obtained from the Grothendieck ring
functor for a group algebra).
The last two of these examples are generalized in Sections 3 and 4
to rings that do not arise from Bouc's construction, and their
definitions appear there.
\end{ex}

The final example in this section is orbifold cohomology of a global quotient.

\begin{ex}\label{orbifold} (Orbifold cohomology.)
Let $L=G$ and let $Y$ be a complex manifold with a left
action of a finite group $G$.
Let $A(g):=\coh^*(Y^g)$, the singular cohomology of the submanifold $Y^g$
of elements invariant under $g$.
The map $c_{g,x}:H^*(Y^g)\rightarrow H^*(Y^{{}^xg})$ 
is induced by the action of $x$ on $Y$.
Let
$$m_{g,h}(\alpha_g,\beta_h)=i_*\left(\alpha_g|_{Y^{\langle g,h\rangle}}
\cdot \beta_h|_{Y^{\langle g,h\rangle}} \cdot c(g,h)\right)$$
where $i:Y^{\langle g,h\rangle}\rightarrow Y^{gh}$ is the natural inclusion,
the pushforward $i_*$ is defined via Poincar\'{e} duality,
and $c(g,h)\in H^*(Y^{\langle g,h\rangle})$ are particular classes
(see \cite{fantechi-gottsche01} or \cite{uribe01} for the details).
By \cite[Theorem 1.18]{fantechi-gottsche01} (see also 
\cite{uribe01}), properties (H1)--(H4) hold,
and $A$ of our Theorem \ref{invariants}(i) is the
associative algebra $\coh^*(Y,G)$ of 
\cite{fantechi-gottsche01} (see also \cite{uribe01}).
The subring $A^G$ is the orbifold cohomology $\coh^*_o([Y/G])$ of 
the orbifold $[Y/G]$.
\end{ex}

The product on orbifold cohomology in the above example
differs from the product on $A^G$ given in our Theorem \ref{invariants}(ii)
by scalar factors (see (\ref{comparison})).
Thus we have two different products on $A^G=\coh_o^*([Y/G])$.
If orbifold cohomology were defined more generally for a quotient of a
topological space by a finite group action,
it may be that there would be
no obvious ring $A$ from which to take the subring $A^G$ of invariants
as the orbifold cohomology.
In addition, if coefficients are taken in a ring $\kf$ of positive
characteristic, the product of Theorem \ref{invariants}(i) may often
be 0 on the invariants $A^G$, due to the (integer) scalar factors that
appear in (\ref{comparison}).
Under these conditions, Theorem \ref{invariants}(ii) may yield the
desired cohomology ring directly.

To illustrate these issues, we consider 
the classifying space $BG=EG/G$ of the finite
group $G$ (see \cite{benson91b} for the definition).
Given $a,b,c\in G$ with $ab=c$, there is an
arrow in $EG$ from $b$ to $c$ labeled by $a$.
A fourth element $g\in G$ acts on the corresponding cell in $EG$ by sending
it to an arrow from $gb$ to $gc$ labeled by $gag^{-1}$.
This action commutes with the action of $G$ on $EG$ by right multiplication
by elements of $G$, and so yields an action of $G$ on $BG$.
The invariant subspace $(BG)^g$ of $g\in G$ may be identified with
$BC(g)$, the classifying space of the centralizer of $g$ in $G$.

Additively, the orbifold cohomology of $BG$ should thus be
$$\coh^*_o([BG/G])=\biggl(\bigoplus_{g\in G}\coh^*(BC(g))\biggr)^G
  \cong\bigoplus_{i=1}^t \coh^*(BC(g_i)),$$
as the group $C(g)$ acts trivially on $\coh^*(BC(g))\cong \coh^*(C(g))$.
(Here $g_1,\ldots,g_t$ is a set of representatives of the conjugacy
classes of $G$.)
If coefficients are taken in a fixed field $\kf$, we may rewrite this as
$$\coh^*_o([BG/G])\cong \bigoplus_{i=1}^t\coh^*(C(g_i))\cong \HH^*(\kf G).$$
That is, $\coh_o^*([BG/G])$ is additively isomorphic to the Hochschild
cohomology of the group algebra $\kf G$.
We may now use the definition of $m_{g,h}$ 
from Example \ref{hh*LG} (with $L=G$), 
making $\coh^*_o([BG/G])$ isomorphic to $\HH^*(\kf G)$ as rings.

\section{Hochschild cohomology of group crossed products}

In this section we prove a product formula (Theorem \ref{productformula})
for the Hochschild cohomology ring of a group crossed product, generalizing
the formula \cite[Thm.\ 5.1]{siegel-witherspoon99} for group algebras.
As a consequence, such a Hochschild cohomology ring satisfies the
properties of Section 2, and its product formula coincides with
that of Theorem \ref{invariants}(ii).

Let $G$ be a finite group acting by automorphisms on an algebra
$S$ over a {\em field} $\kf$.\footnote{Working over a field simplifies
Hochschild cohomology.}
Let $\sigma : G\times G\rightarrow Z(S)^{\times}$ (the units in the
center of $S$) be a two-cocycle, that is
\begin{equation}\label{cocycle}
{}^g\sigma(h,k)\sigma(g,hk)=\sigma(g,h)\sigma(gh,k)
\end{equation}
for all $g,h,k\in G$.\footnote{We require the image of $\sigma$ 
to be central in $S$ to be consistent with having an action of $G$ on $S$.
More generally, the action of $G$ could be twisted by $\sigma$
with noncentral image.
See for example \cite[Chapter 7]{montgomery93}.}
Assume that $\sigma$ is normalized, that is
$\sigma(g,1)=1=\sigma(1,g)$ for all $g\in G$.

The {\em crossed product algebra} $R:=S\#_{\sigma}G$ (or
$S\#_{\sigma} \kf G$) has underlying
vector space $S\otimes \kf G$ (where $\otimes=\otimes_{\kf}$)
and multiplication given by
$$(s\otimes g)(t\otimes h)=s({}^g t)\sigma(g,h)\otimes gh$$
for all $s,t\in S$ and $g,h\in G$.
We will write $s\overline{g}:=s\otimes g$ to shorten notation.
Note that $S$ is a subalgebra of $R$, but $\kf G$ is not.
For each $g\in G$,
the element $\overline{g}\in R$ is invertible with
$$(\overline{g})^{-1}=\sigma^{-1}(g,g^{-1})\overline{g^{-1}}.$$
The action of $G$ on $S$ becomes an inner action on $R$, as
\begin{equation}\label{inner}
{}^g s = \overline{g} s (\overline{g})^{-1}
\end{equation}
for all $g\in G, s\in S$.

A number of authors have studied Hochschild (co)homology of such
group crossed products or more general algebras.
Sanada \cite{sanada93} and \c{S}tefan \cite{stefan95} gave spectral
sequences relating Hochschild cohomology to that of the components 
$S$ and $G$.
(\c{S}tefan treated the more general case of a Hopf Galois extension.)
Cornick \cite{cornick95} and Lorenz \cite{lorenz92} gave spectral
sequences to describe components of Hochschild homology of group-graded
algebras that are indexed by the conjugacy classes of $G$.
Here we will describe the overall product structure
of the Hochschild cohomology ring $\HH^*(R)$ {\em in terms of} components
indexed by conjugacy classes.

We will first give a decomposition of the Hochschild cohomology ring
$\HH^*(R)$ into such components.
We recall the definition of Hochschild cohomology: 
As $R$ is an algebra over a field, its Hochschild cohomology may be
defined as
$$\HH^*(R):=\Ext^*_{R^e}(R,R),$$
where $R^e:=R\otimes R^{op}$ acts on the left on $R$ by left and right
multiplication.
More generally, the Hochschild cohomology of $R$ with coefficients in an 
$R$-bimodule $M$ is $\HH^*(R,M):=\Ext^*_{R^e}(R,M)$.
(So $\HH^*(R)=\HH^*(R,R)$.)
This may be expressed in terms of the (acyclic) Hochschild complex:
\begin{equation}\label{acyclic}
\cdots\stackrel{\delta_3}{\longrightarrow} R^{\otimes 4}
\stackrel{\delta_2}{\longrightarrow} R^{\otimes 3}\stackrel{\delta_1}
{\longrightarrow} R^e \stackrel{m}{\longrightarrow} R\longrightarrow 0
\end{equation}
is an $R^e$-free resolution of $R$, where $m$ is the multiplication map and
$$\delta_n(r_0\otimes r_1\otimes \cdots\otimes r_{n+1})=
\sum_{i=0}^n(-1)^ir_0\otimes \cdots\otimes r_ir_{i+1}\otimes\cdots\otimes
r_{n+1}.$$
Dropping the term $R$ from the complex (\ref{acyclic}) above, and applying
$\Hom_{R^e}(-,M)$, we have the Hochschild (cochain) complex
\begin{equation}\label{cochain}
\Hom_{R^e}(R^e,M)\stackrel{\delta_1^*}{\longrightarrow}
\Hom_{R^e}(R^{\otimes 3},M)\stackrel{\delta_2^*}{\longrightarrow}
\Hom_{R^e}(R^{\otimes 4},M)\stackrel{\delta_3^*}{\longrightarrow}\cdots
\end{equation}
Thus $\HH^n(R,M)=\Ker(\delta_{n+1}^*)/\Ima(\delta_n^*)$.

If $H$ is a subgroup of $G$, let $$\Delta(H):=\bigoplus_{h\in H}S\overline{h}
\otimes S(\overline{h})^{-1},$$ a subalgebra of $R^e$, and let
$\Delta:=\Delta(G)$.
Let $g_1,\ldots,g_t$ be a set of representatives of the 
conjugacy classes of $G$.

\begin{lem}\label{hh*-additive}
Let $R=S\#_{\sigma}G$.
There is an isomorphism of vector spaces 
$$\HH^*(R)\cong \bigoplus_{i=1}^t\Ext^*_{\Delta
  (C(g_i))}(S,S\overline{g_i}).$$
\end{lem}

\begin{proof}
Note that $R\cong R^e\otimes_{\Delta}S$ as an
$R^e$-module (see \cite[Lemma 3.3]{boisen92}), 
and so by the Eckmann-Shapiro Lemma \cite[Cor.\ 2.8.4]{benson91} we have
$$\HH^*(R)\cong \Ext^*_{R^e}(R^e\otimes_{\Delta}S,R)\cong \Ext
  ^*_{\Delta}(S,R).$$
Now, as a $\Delta$-module, $R\cong\oplus_{i=1}^tR_{{\mathcal C}_i}$,
where ${\mathcal C}_i$ is the conjugacy class of $g_i$ and $R_{{\mathcal C}_i}
=\oplus_{h\in {\mathcal C}_i}S\overline{h}$.
Therefore $\HH^*(R)\cong\oplus_{i=1}^t\Ext^*_{\Delta}(S,R_{{
\mathcal C}_i})$. 

As a $\Delta$-module, $R_{{\mathcal C}_i}$ is isomorphic to 
the coinduced module $\Hom_{\Delta(C(g_i))}(\Delta,S\overline{g_i})$, where 
the action of $\Delta$ is given by $((r_0\otimes r_1)\cdot f)(s_0\otimes s_1) 
=f(s_0r_0\otimes r_1s_1)$ for $f\in\Hom_{\Delta(C(g_i))}(\Delta,S\overline
{g_i})$:
The map $R_{{\mathcal C}_i}\rightarrow \Hom_{\Delta (C(g_i))}(\Delta, 
S\overline{g_i})$ defined by
$$(\overline{h})r(\overline{g})(\overline{h})^{-1}\mapsto
\left(((\overline{h})^{-1}\otimes\overline{h})\mapsto 
r\overline{g}\right)$$
is an isomorphism of $\Delta$-modules, where $h$ ranges over a set of coset
representatives for $C(g_i)$ in $G$.
Another way to view this isomorphism is to notice that $R_{{\mathcal {C}}_i}
\cong \Delta\otimes_{\Delta(C(g_i))}S\overline{g_i}$, and this induced
module is isomorphic to the coinduced module $\Hom_{\Delta(C(g_i))}
(\Delta,S\overline{g_i})$.
Passman has provided us with a proof of the fact that the induced
and coinduced modules for finite group crossed products are isomorphic.

Again by the Eckmann-Shapiro Lemma,
$$\Ext^*_{\Delta}(S,R_{{\mathcal C}_i})\cong \Ext^*_{\Delta(C(g_i))}(S,
S\overline{g_i}),$$ 
which proves the lemma.
\end{proof}

We will show that $\HH^*(R)$ satisfies the properties of
Section 2, and in
the process will determine how the Hochschild cup product behaves with
respect to the additive decomposition of Lemma \ref{hh*-additive}.
We will need to develop some general theory regarding such Hochschild 
cohomology rings first.

Let $M$ be a $\Delta(H)$-module for a subgroup $H\leq G$, and $K$ a
subgroup of $H$.
Let 
\begin{equation}\label{PS}
\cdots \rightarrow P_2\rightarrow P_1\rightarrow P_0\rightarrow S\rightarrow 0
\end{equation}
be a $\Delta(H)$-projective resolution of $S$.
Since $\Delta(H)$ is free as a $\Delta(K)$-module, this restricts to a
$\Delta(K)$-projective resolution of $S$, and so there is a
{\em restriction} map
$$\res^H_K:\Ext^*_{\Delta(H)}(S,M)\rightarrow \Ext^*_{\Delta(K)}(S,M),$$
defined at the chain level by the inclusion 
$\Hom_{\Delta(H)}(P_n,M)\rightarrow
\Hom_{\Delta(K)}(P_n,M)$. 
If $g\in G$, conjugation by $(\overline{g})^{-1}$ 
is an algebra isomorphism from $\Delta( {}^gH)$ to $\Delta(H)$.
Denote by ${}^gM$ the $\Delta({}^gH)$-module corresponding to $M$
under this algebra isomorphism. 
Then we have {\em conjugation} maps
$$g^*:\Ext^*_{\Delta(H)}(S,M)\rightarrow \Ext^*_{\Delta({}^gH)}(S,{}^gM)$$
induced by the algebra isomorphism $\Delta({}^gH)\stackrel{\sim}
{\longrightarrow} \Delta(H)$.
They may be described
at the chain level as follows.  Let 
$$\cdots\rightarrow {}^gP_2\rightarrow {}^gP_1\rightarrow
{}^gP_0\rightarrow S\rightarrow 0$$
be the $\Delta({}^gH)$-projective resolution of $S$ corresponding to
(\ref{PS}).
Denote by ${}^gp$ an element of ${}^gP_n$ thus corresponding to $p\in P_n$.
We let $g^*(f)({}^gp)={}^gf(p)$.
Note that $h\in H$ acts trivially by conjugation on $\Ext^*_{\Delta(H)}
(S,M)$.  The {\em corestriction} map
$$\cores^H_K:\Ext^*_{\Delta(K)}(S,M)\rightarrow \Ext^*_{\Delta(H)}(S,M)$$
is defined at the chain level
by $\cores^H_K(f)(p)=\sum_{g\in H/K} gf(g^{-1}p)$.
Strictly speaking, this sum is over a set of coset representatives of
$H/K$, but since $K$ acts trivially on $\Hom_{\Delta(K)}(P_n,M)$, 
the choice of representatives does not matter.
This induces a map on cohomology as it commutes with the cochain
maps (see for example \cite[Section 3.6]{benson91b} for the case of
group cohomology).

We recall the Hochschild cup product on $\HH^*(R)$,
defined at the chain level on the Hochschild complex (\ref{cochain})
when $M$ is itself a ring (see \cite{gerstenhaber63}):
If $f\in \Hom_{R^e}(R^{\otimes m},M)$ and $f'\in\Hom_{R^e}(R^{\otimes n},M)$,
then $f\otimes f'\in\Hom_{R^e}(R^{\otimes (m+n)},M)$ is defined by
$$(f\otimes f')(r_1\otimes\cdots\otimes r_{m+n})=f(r_1\otimes \cdots
\otimes r_m) f'(r_{m+1}\otimes\cdots\otimes r_{m+n}).$$
In order to translate this into a product on $\Ext^*_{\Delta}(S,R)$, we 
will give a $\Delta$-projective resolution of $S$ that induces to the
Hochschild complex (\ref{acyclic}) for $R$:
For each $n\geq 0$,
let $\Delta_n$ be the $\Delta$-submodule of $R^{\otimes (n+2)}$
consisting of sums of all $s_0\overline{g_0}\otimes\cdots\otimes s_{n+1}
\overline{g_{n+1}}$ ($s_i\in S$, $g_i\in G$) such that 
$g_0\cdots g_{n+1}=1$ in $G$.
Thus $\Delta_0=\Delta$, and each $\Delta_n$ is a projective $\Delta$-module.
Then
\begin{equation}\label{subcomplex}
\cdots\stackrel{\delta_3}{\longrightarrow}\Delta_2\stackrel{\delta_2}
{\longrightarrow}\Delta_1\stackrel{\delta_1}{\longrightarrow}\Delta_0
\stackrel{m}{\longrightarrow} S\longrightarrow 0
\end{equation}
is a projective resolution of the $\Delta$-module $S$, where the maps
are restrictions of the maps from the Hochschild complex (\ref{acyclic}).
The above complex (\ref{subcomplex}) is exact, as there is a chain contraction
$s_n:\Delta_{n-1}\rightarrow\Delta_n$ given by 
$$s_n(r_0\otimes\cdots\otimes
r_n)=r_0\otimes\cdots\otimes r_n\otimes 1.$$
The Hochschild complex (\ref{acyclic}) is just the complex (\ref{subcomplex})
induced from $\Delta$ to $R^e$.
Thus the isomorphism $\Ext_{R^e}^n(R,R)\stackrel{\sim}{\rightarrow}
\Ext^n_{\Delta}(S,R)$ of the Eckmann-Shapiro Lemma
is given at the chain level simply by
restricting maps from $\Hom_{R^e}(R^{\otimes (n+2)},R)$ to
$\Hom_{\Delta}(\Delta_n,R)$.
The Hochschild cup product on $\HH^*(R)$ now yields a cup product on
$\Ext^*_{\Delta}(S,R)$ as follows.
If $f\in\Hom_{\Delta}(\Delta_m,R)$ and $f'\in\Hom_{\Delta}(\Delta_n,R)$
then $f\otimes f'\in\Hom_{\Delta}(\Delta_{m+n},R)$ is defined by
$$(f\otimes f')(r_1\otimes\cdots\otimes r_{m+n})=\overline{g} f(r_1'\otimes
r_2\otimes\cdots\otimes r_m)f'(r_{m+1}\otimes\cdots \otimes
r_{m+n-1}\otimes r_{m+n}') (\overline{g})^{-1},$$
where $r_1=\overline{g}r_1'$ and $r_{m+n}=r'_{m+n}(\overline{g})^{-1}$ are
chosen so that $r_1'\otimes r_2\otimes\cdots\otimes r_m$ and
$r_{m+1}\otimes\cdots\otimes r_{m+n-1}\otimes r_{m+n}'$ are in $\Delta_m$
and $\Delta_n$, respectively.
Similarly, we have a Hochschild cup product on $\Ext^*_{\Delta(H)}(S,R)$
for any subgroup $H$ of $G$.

\begin{lem}\label{FM} 
Let $K$ and $H$ be subgroups of $G$.
The following relations hold for all $\alpha\in\Ext^*_{\Delta(H)}
(S,R)$ and $\beta\in\Ext^*_{\Delta(K)}(S,R)$:
\begin{itemize}
\item[(i)] (Frobenius property) If $K<H$ then 
$$\cores^H_K(\res^H_K(\alpha)\smile\beta)=\alpha\smile\cores^H_K(\beta) \
\mbox{ and } \ \cores^H_K(\beta\smile\res^H_K(\alpha))=\cores^H_K(\beta)
\smile\alpha.$$
\item[(ii)] (Mackey property)
$$\res^G_K(\cores^G_H(\alpha))=\sum_{x\in D}\cores^K_{K\cap {}^xH}
(\res^{{}^xH}_{K\cap {}^xH}(x^*\alpha)),$$
where $D$ is a set of double coset representatives for $K\backslash
G/H$.
\end{itemize}
\end{lem}

\begin{proof}
Property (i) holds at the chain level by 
the definitions, and (ii) follows from a particular ordering of $G/H$ (see 
\cite[Thm.\ 4.2.6]{evens91}).
\end{proof}

In fact the assignment of $\Ext^*_{\Delta(H)}(S,M)$ to each subgroup $H$
of $G$ is a Green functor (see \cite{thevenaz??} for the definition).
The additional properties that we will need are straightforward: 
\begin{eqnarray}
&&g^*\circ \res^H_K=\res^{{}^gH}_{{}^gK}\circ g^* \mbox{ and }
g^*\circ\cores^H_K=\cores^{{}^gH}_{{}^gK}\circ g^* \
(g\in G \mbox{ and } K<H<G), \mbox{ and}\\
&&\hspace{.2cm}\res^K_L\circ \res^H_K =\res^H_L
\mbox{ and }\cores^H_K\circ\cores^K_L=\cores^H_L\mbox{ whenever } L<K<H<G.
\end{eqnarray}

We will need explicit maps giving the isomorphism of Lemma \ref{hh*-additive}.
For each $g\in G$, the inclusion map $\theta_g:S\overline{g}\rightarrow R$
and projection map $\pi_g:R\rightarrow S\overline{g}$ are both $\Delta
(C(g))$-homomorphisms.
Thus they induce the following maps on cohomology:
\begin{eqnarray*}
\theta_g^* : &:& \Ext^*_{\Delta(C(g))}(S,S\overline{g})\rightarrow
  \Ext^*_{\Delta(C(g))}(S,R)\\
\mbox{and } \ \pi_g^* &:& \Ext^*_{\Delta(C(g))}(S,R)\rightarrow
   \Ext^*_{\Delta(C(g))}(S,S\overline{g}).
\end{eqnarray*}
Let $g_1,\ldots,g_t$ be a set of representatives of conjugacy classes of
$G$, and write $\theta_i:=\theta_{g_i}$ and $\pi_i:=\pi_{g_i}$.
The isomorphism $\Ext^*_{\Delta}(S,R)\stackrel{\sim}{\rightarrow}
\oplus_{i=1}^t\Ext^*_{\Delta(C(g_i))}(S,S\overline{g_i})$ 
is given explicitly by
$\zeta\mapsto (\pi_i^*\res^G_{C(g_i)}(\zeta))_i$, and
its inverse sends $\alpha\in\Ext^*_{\Delta(C(g_i))}(S,S\overline{g_i})$ to 
$\cores^G_{C(g_i)}\circ\theta_i^*(\alpha)$.
(See \cite[Lemmas 4.1 and 4.2]{siegel-witherspoon99} for more details
in the case of a group algebra.)

Let $g,h\in G$.  The following further properties are straightforward:
\begin{multline}
\mbox{If }W\leq C(g),\mbox{ then }h^*\circ\theta_g^*=\theta_{{}^hg}^*\circ
h^*\mbox{ and }h^*\circ\pi_g^*=\pi_{{}^hg}^*\circ h^*\\
\mbox{ as maps from }\Ext^*_{\Delta(W)}(S,S\overline{g})\mbox{ to }\Ext^*
_{\Delta({}^hW)}(S,R).
\end{multline}
\begin{equation}
\mbox{If }W'\leq W\leq C(g),\mbox{ then }\theta_g^*\mbox{ and }
\pi_g^*\mbox{ commute with }\res^W_{W'}\mbox{ and }\cores^W_{W'}.\hspace{.45in}
\end{equation}
\begin{multline}
\mbox{If }W\leq C(g)\cap C(h),\mbox{ then }\pi_g^*\circ\theta_h^*=
\delta_{g,h}\id\\
\mbox{ as maps from }\Ext^*_{\Delta(W)}(S,S\overline{h})
\mbox{ to }\Ext^*_{\Delta(W)}(S,S\overline{g}).
\end{multline}

We now have the following product formula. 
For $\alpha\in\Ext^*_{\Delta(C(g_i))}(S,S\overline{g_i})$,
write $\gamma_i(\alpha)=\cores^G_{C(g_i)}\circ\theta_{i}^*(\alpha)$,
which is the image of $\alpha$ in
$\Ext^*_{\Delta}(S,R)\cong \HH^*(R)$ under the isomorphism of Lemma
\ref{hh*-additive}.

\begin{thm}\label{productformula}
Let $\alpha\in\Ext^*_{\Delta(C(g_i))}(S,S\overline{g_i}),
\beta\in\Ext^*_{\Delta(C(g_j))}(S,S\overline{g_j})$, and
$\gamma_i(\alpha)$ and $\gamma_j(\beta)$ their images in 
$\Ext_{\Delta}^*(S,R)\cong \HH^*(R)$ under the isomorphism of Lemma
\ref{hh*-additive}.
Then the cup product of $\gamma_i(\alpha)$ and $\gamma_j(\beta)$ in
$\HH^*(R)$ is given by
$$\sum_{x\in D}
\gamma_k\left(\cores^{C(g_k)}_{{}^yC(g_i)\cap {}^{yx}C(g_j)}\pi_{k}^*\left(
  \theta_{{}^yg_i}^*\res^{{}^yC(g_i)}_{{}^yC(g_i)
 \cap{}^{yx}C(g_j)}y^*\alpha \smile \theta^*_{{}^{yx}g_j}
\res^{{}^{yx}C(g_j)}_{{}^yC(g_i)\cap{}^{yx}C(g_j)}(yx)^*\beta\right)\right),$$
where $D$ is a set of representatives of double cosets $C(g_i)\backslash
G/C(g_j)$, and $k=k(x)$ and $y=y(x)$ are chosen so that $g_k={}^yg_i
{}^{yx}g_j$.
\end{thm}

\begin{proof}
By the Frobenius and Mackey properties (Lemma \ref{FM}), 
$\gamma_i(\alpha)\smile\gamma_j(\beta)$ is equal to 
\begin{eqnarray*}
&& \hspace{-1in}\cores^G_{C(g_i)}\left(\theta^*
  _{i}\alpha\right)\smile \cores^G_{C(g_j)}\left(\theta_{j}^*\beta\right)\\
  &=& \cores^G_{C(g_i)}\left(\theta^*_{i}\alpha\smile\res^G_{C(g_i)}\cores^G
  _{C(g_j)}\theta^*_{j}\beta\right)\\
  &=&\sum_{x\in D}\cores^G_{C(g_i)}\left(\theta^*_{i}\alpha\smile
  \cores^{C(g_i)}_{C(g_i)\cap {}^xC(g_j)}\res^{{}^xC(g_j)}_{C(g_i)\cap
  {}^xC(g_j)} x^*\theta^*_{j}\beta\right)\\
  &=&\sum_{x\in D}\cores^G_{C(g_i)\cap {}^xC(g_j)}
  \left(\res^{C(g_i)}_{C(g_i)\cap {}^xC(g_j)}
  \theta^*_{i}\alpha\smile\res^{{}^xC(g_j)}_{C(g_i)\cap {}^xC(g_j)}
 x^*\theta^*_{j}\beta\right)
\end{eqnarray*}
Inserting the identity map $\id=\sum_{k=1}^t\gamma_k\pi_{k}^*\res^G_{C(g_k)}$,
and applying the Mackey property and properties (3.9)--(3.13), 
$\gamma_i(\alpha)\smile\gamma_j(\beta)$ equals
\begin{eqnarray*}
&&\hspace{-.4in} 
  \sum_{k=1}^t\sum_{x\in D}\gamma_k
  \pi^*_{k}\res^G_{C(g_k)}\cores^G_{C(g_i)\cap {}^xC(g_j)}\left(\res
  ^{C(g_i)}_{C(g_i)\cap{}^xC(g_j)}\theta^*_{i}\alpha\smile \res
  ^{{}^xC(g_j)}_{C(g_i)\cap {}^xC(g_j)}x^*\theta^*_{j}\beta\right)\\
 &=&\sum_{k=1}^t\sum_{x,y}\gamma_k\left(\pi^*_{k}\cores^{C(g_k)}_{W}
  \res_{W}^{{}^yC(g_i)\cap {}^{yx}C(g_j)}y^*\left(\res^{C(g_i)}_{C(g_i)\cap
   {}^xC(g_j)}\theta^*_{i}\alpha\smile\res^{{}^xC(g_j)}_{C(g_i)\cap 
   {}^xC(g_j)}x^*\theta^*_{j}\beta\right)\right)\\
  &=& \sum_{k=1}^t\sum_{x,y}\gamma_k\left(\cores^{C(g_k)}_{W}\pi^*_{k}\left(
  \theta^*_{{}^yg_i}\res^{{}^yC(g_i)}_{W} y^*\alpha\smile\theta^*
  _{{}^{yx}g_j}\res^{{}^{yx}C(g_j)}_{W}(yx)^*\beta\right)\right),
\end{eqnarray*}
where $y$ runs over a set of double coset representatives for 
$C(g_k)\backslash G/C(g_i)\cap {}^xC(g_j)$, and 
$W=C(g_k)\cap {}^yC(g_i)\cap {}^{yx}C(g_j)$.
Now $\theta^*_{{}^yg_i}\res^{{}^yC(g_i)}_{W}
y^*\alpha\smile\theta^*_{{}^{yx}g_j}\res^{{}^{yx}C(g_j)}_{W} (yx)^*\beta$
is in the image of the map $\theta^*_{{}^yg_i{}^{yx}g_j}$ from
$\Ext^*_{\Delta(W)}(S,S\overline{{}^yg_i {}^{yx}g_j})$ to 
$\Ext^*_{\Delta(W)}(S,R)$, and so if we apply $\pi_{g_k}^*$, this can
only be nonzero when $g_k={}^yg_i{}^{yx}g_j$.
But each $x$ determines a unique $k$ and double coset $C(g_k)y(C(g_i)\cap
{}^xC(g_j))$ for which this holds.
So we may take $k=k(x)$, $y=y(x)$, and then 
${}^yC(g_i)\cap {}^{yx}C(g_j)\leq C(g_k)$,
so $W= {}^yC(g_i)\cap {}^{yx}C(g_j)$.
\end{proof}

Finally we show how $\HH^*(R)$ fits into the construction of Section 2.
For each $g\in G$, let $A(g):=\Ext_{\Delta(C(g))}^*(S,S\overline{g})$.
The conjugation maps $c_h=h^*$ 
yield isomorphisms of $k$-modules $A(g)\stackrel{\sim}
{\rightarrow}A({}^hg)$.  Define
$m_{g,h}:A(g)\times A(h)\rightarrow A(gh)$ by
$$m_{g,h}(\alpha_g,\beta_h)=\cores^{C(gh)}_{C(g)\cap C(h)}
\left(\pi^*_{gh}\left(\theta^*_g\res^{C(g)}_{C(g)\cap C(h)}
\alpha_g\smile \theta^*_h
\res^{C(h)}_{C(g)\cap C(h)}\beta_h\right)\right).$$
Then properties (H1)--(H3) of Section 2 are straightforward.
A comparison of Corollary \ref{product} and Theorem \ref{productformula}
shows that property (H4$'$) holds, so that
Theorem \ref{invariants}(ii) yields a product on $A^G$.
By (\ref{decomposition}), Corollary \ref{product}, Lemma \ref{hh*-additive}
and Theorem \ref{productformula}, there is an algebra isomorphism
$$A^G\cong \bigoplus_{i=1}^t \Ext_{\Delta(C(g_i))}(S,S\overline{g_i})
\cong\HH^*(R).$$

\section{Representations of abelian extensions of Hopf algebras}

In this section we prove a product formula (Theorem \ref{tensor}) for
the Grothendieck ring of modules for an abelian extension of Hopf algebras.
As a consequence, this Grothendieck ring satisfies the properties of
Section 2, and its product formula coincides with that of Theorem 
\ref{invariants}(ii).
This generalizes a known formula
for modules of the Drinfel'd double of a finite group (or equivalently
of Hopf bimodules of the group algebra \cite{cibils97}).
Related results are in \cite{cibils97,witherspoon96}.

We assume in this section that $\kf$ is an {\em algebraically closed
field} as we will use Schur's Lemma, 
but we do not place any restrictions on the characteristic of $\kf$.
All our modules will be finite dimensional.

Let $L$ be a finite group acting on another finite group $G$.  
This induces an action of $L$ on the linear dual $(\kf G)^*$
of the group algebra $\kf G$.
For each $g\in G$, let $p_g\in (\kf G)^*$ denote the function dual to
$g$ in the basis $G$ of $\kf G$, that is
$$p_g(h):=\delta_{g,h} \ \ (h\in G).$$
We have ${}^xp_g=p_{{}^xg}$ for all $x\in L, g\in G$.
Let $\sigma : L\times L\rightarrow ((\kf G)^*)^{\times}$ and 
$\tau:L\rightarrow ((\kf G)^*)^{\times}\times 
((\kf G)^*)^{\times}$ be maps giving rise to a (quasi) Hopf algebra
$H:=(\kf G)^*\#_{\sigma}^{\tau}L$ (or $(\kf G)^*\#_{\sigma}^{\tau}\kf L$)
as follows.
As an algebra, $H$ is the crossed product $(\kf G)^*\#_{\sigma}L$
defined in Section 3.
Again we write $p_g\overline{x}:=p_g\otimes x$ to shorten the notation;
$(\kf G)^*$ will be a Hopf subalgebra of $H$, whereas $\kf L$ will not.
For each $x,y\in L$, let
\begin{equation}\label{st}
\sigma(x,y)=\sum_{g\in G}\sigma_g(x,y)p_g, \ \mbox{ and } \
\tau(x)=\sum_{g,h\in G}\tau_{g,h}(x)p_g\otimes p_h,
\end{equation}
for scalars $\sigma_g(x,y)$ and $\tau_{g,h}(x)$.
Then the product in $H$ may be written
$$(p_g\overline{x})\cdot(p_h\overline{y})=
 \delta_{g, {}^xh} \sigma_{g}(x,y)p_g\overline{xy}.$$
Define the coproduct by\footnote{This is the case that $H$ is a
{\em cocentral} abelian extension. In general, the coproduct is more
complicated.}
$$\Delta(p_g\overline{x}):=\sum_{\substack{h,k\in G\\hk=g}}
\tau_{h,k}(x)p_h\overline{x}\otimes p_{k}\overline{x}.$$
In order that $H=(\kf G)^*\#^{\tau}_{\sigma}L$ be a (quasi) Hopf
algebra, $\tau$ (as well as $\sigma$) must satisfy certain properties.
See \cite{andruskiewitsch96} for the general case of a Hopf algebra,
\cite{kashina-mason-montgomery02} for (cocentral)
abelian extensions in particular,
\cite{chari-pressley94} for general facts about quasi Hopf algebras,
and the references given at the end of the following example for a
special case.

\begin{ex} 
Let $G$ act on itself by conjugation.
Let $\omega: G\times G\times G\rightarrow \kf^{\times}$ be any 
{\em three}-cocycle,
that is
$$\omega(a,b,c)\omega(a,bc,d)\omega(b,c,d)=\omega(ab,c,d)\omega(a,b,cd)$$
for all $a,b,c,d\in G$.
Assume that $\omega$ is normalized so that $\omega(a,b,c)$ is equal to 1
whenever one of $a$, $b$, or $c$ is 1.
There is an associated {\em two}-cocycle $\sigma:G\times G\rightarrow
(\kf G)^*$ given by (\ref{st}) where
$$
\sigma_g(x,y)=\frac{\omega(g,x,y)\omega(x,y,(xy)^{-1}gxy)}{\omega
(x,x^{-1}gx,y)},$$
and an associated two-cycle $\tau: G\rightarrow (\kf G)^*\times (\kf G)^*$ 
given by (\ref{st}) where
$$\tau_{g,h}(x)=\frac{\omega(g,h,x)\omega(x,x^{-1}gx,x^{-1}hx)}
{\omega(g,x,x^{-1}hx)}.$$
The {\em twisted Drinfel'd} (or {\em quantum}) {\em double} 
$D^{\omega}(G):=(\kf G)^*\#^{\tau}_{\sigma} G$ is a quasi Hopf algebra.
In case $\omega$ is a coboundary, this is isomorphic to the
{\em Drinfel'd double} $D(G):= (\kf G)^*\# G$, a Hopf algebra.
These (quasi) Hopf algebras and their representations appear in 
\cite{bantay91,dijkgraaf-pasquier-roche90,mason??,witherspoon96,witherspoon??}.
\end{ex} 

We will first apply Clifford theory to the crossed products 
$H=(\kf G)^*\#^{\tau}_{\sigma} L$ to obtain a description 
of all simple $H$-modules.
Such a description was first found in \cite{kashina-mason-montgomery02}
by more direct methods in case $H$ is semisimple.
Note that we do not need the coalgebra structure of $H$ to
describe $H$-modules.
Clifford theory for group-graded rings generally
is developed in \cite{dade70,dade86}, and
more specifically for crossed products the theory is in 
\cite[\S11C]{curtis-reiner__}. However we will use the terminology and
notation of \cite{montgomery-witherspoon98}, recalling the needed
results as we go.  

Up to isomorphism, the simple $(\kf G)^*$-modules are the ideals 
$\kf p_g$ $(g\in G)$.
The {\em stabilizer}
of the $(\kf G)^*$-module $\kf p_g$ in $L$ is the subgroup of all $x\in L$
such that ${}^x(\kf p_g)\cong \kf p_g$. This is then the subgroup
$$L_g=\{x\in L\mid {}^x g=g\}.$$
By (\ref{cocycle}), the function 
$\sigma_g:G\times G\rightarrow \kf^{\times}$ defined by (\ref{st})
restricts to a two-cocycle $L_g\times L_g\rightarrow\kf^{\times}$.
Therefore we may form the {\em twisted group algebra} 
$$\kf_{\sigma_g}L_g :=
\kf \#_{\sigma_g} L_g,$$ 
that is the algebra with basis
$\{\overline{x}\mid x\in L_g\}$ and multiplication $\overline{x}\cdot
\overline{y}=\sigma_g(x,y)\overline{xy}$.
This is in fact isomorphic to the subalgebra $p_g\overline{L_g}$
of $H$ generated by
all $p_g\overline{x}$ ($x\in L_g$), so no confusion should arise
from the choice of notation.

For each $g\in G$, we define the subalgebra of $H$,
$$H_g:=(\kf G)^*\#_{\sigma} L_g.$$
Let $H_g\otimes_{(\kf G)^*}\kf p_g$ be the $H_g$-module induced from
the $(\kf G)^*$-module $\kf p_g$.
Let
$$E:=\End_{H_g}(H_g\otimes_{(\kf G)^*}\kf p_g)^{op},$$
the endomorphism algebra of the module $H_g\otimes_{(\kf G)^*}\kf p_g$,
taken with multiplication opposite that of 
composition of functions.\footnote{We must take the opposite multiplication
as we write our functions on the left.}

\begin{lem}\label{endo}
There is an algebra isomorphism
$E\cong \kf_{\sigma_g}L_g$.
\end{lem}

\begin{proof} There are isomorphisms of vector spaces
$$\End_{H_g}(H_g\otimes_{(\kf G)^*}\kf p_g)\cong \Hom_{(\kf G)^*}
(\kf p_g,H_g\otimes
_{(\kf G)^*}\kf p_g)\cong\bigoplus_{x\in L_g}\Hom_{(\kf G)^*}(\kf p_g,
\overline{x}\otimes \kf p_g).$$
Each such endomorphism is thus determined by the image of the element
$p_g=1\otimes p_g$, 
and this must be a linear combination of the $\overline{x}
\otimes p_g$ $(x\in L_g)$.  The product of the functions $\phi_{x},
\phi_{y}$, where $\phi_{x}(p_g)=\overline{x}\otimes p_g$
and $\phi_{y}(p_g)=\overline{y}\otimes p_g$, is given
by the opposite of composition:
$$(\phi_{y}\circ \phi_{x})(p_g)=
\phi_{y}(\overline{x}\otimes p_g)=\overline{x}\phi_{y}
(1\otimes p_g)=\overline{x}(\overline{y}\otimes p_g)=\sigma
(x,y)\overline{xy}\otimes p_g.$$
As $\sigma(x,y)=\sum_{h\in G}\sigma_h(x,y)p_h$, and the tensor product
is over $(\kf G)^*$, this is in fact equal to $\sigma_g(x,y)\overline{xy}
\otimes p_g$.
Therefore $(\phi_y\circ \phi_x)(p_g)=\sigma_g(x,y)\phi_{xy}(p_g)$.
\end{proof}

The endomorphism algebra $E\cong \kf_{\sigma_g}L_g$ plays a crucial role
in the Clifford correspondence.
(See for example \cite{witherspoon99,witherspoon02}.)
In order to use the results of \cite{montgomery-witherspoon98},
where the algebra $E$ is not explicitly mentioned, we point out that
$$\sigma_g(x,y)=\sigma_g^{-1}(y^{-1},x^{-1})
\sigma_g^{-1}(x^{-1},x)\sigma_g^{-1}(y^{-1},y)
\sigma_g(y^{-1}x^{-1},xy),$$
as follows from (\ref{cocycle}).
Letting $\alpha(x,y)=\sigma_g(y^{-1},x^{-1})$, this shows that $\sigma_g$
is cohomologous to $\alpha^{-1}$.
It may be checked that $\alpha$ is the cocycle of 
\cite[Prop.\ 1.2]{montgomery-witherspoon98} that was used for the
Clifford correspondence there.
As $\sigma_g$ and $\alpha^{-1}$ are cohomologous, there is an algebra
isomorphism $\kf_{\sigma_g}L_g\stackrel{\sim}{\longrightarrow}
\kf_{\alpha^{-1}}L_g$, 
given in this case by $\overline{x}\mapsto \sigma_g^{-1}(x^{-1},x)
\overline{x}$.

Notice that any $H$-module must contain at least one of the simple
$(\kf G)^*$-modules $\kf p_g$, on restriction to $(\kf G)^*$.
Given such an $H$-module, its $H$-submodule generated by $p_g$ is a
direct sum of copies of $\kf p_{{}^xg}$ ($x\in L$) on restriction to
$(\kf G)^*$, as $p_h (\overline{x}p_g)=\delta_{h,{}^xg}p_h\overline{x}$.
Therefore any {\em simple} $H$-module is a direct sum of copies of
conjugates of some $\kf p_g$, on restriction to $(\kf G)^*$.
For the Clifford correspondence, we will also consider $H_g$-modules:
Note that any simple $H_g$-module containing 
$\kf p_g$ on restriction to $(\kf G)^*$
is in fact a direct sum of copies of $\kf p_g$, as $L_g$ is the stabilizer
of $\kf p_g$.

The Clifford correspondence comes in two steps.  In the first step,
there is a bijection between the set of (isomorphism classes of)
simple $\kf_{\sigma_g}L_g$-modules (that is $E$-modules),
and the simple $H_g$-modules
whose restriction to $(\kf G)^*$ contains $\kf p_g$ 
\cite[Prop.\ 1.2]{montgomery-witherspoon98}.  
The second step of the Clifford correspondence
is a bijection between the set of simple $H_g$-modules whose restriction to 
$(\kf G)^*$ contains $\kf p_g$, and the simple
$H$-modules whose restriction to $(\kf G)^*$ contains
$\kf p_g$ \cite[Prop.\ 1.1]{montgomery-witherspoon98}.
This second step is simply given by tensor induction of modules.
An explicit description of such $H$-modules is given in
\cite[Thm.\ 1.3]{montgomery-witherspoon98}. In our situation, it
provides a new proof of the following proposition which appears
as \cite[Theorem 3.3]{kashina-mason-montgomery02} in case 
$H$ is semisimple.
The action described in the proposition below looks
simpler than that of \cite[Theorem 1.3]{montgomery-witherspoon98}
due to application of the algebra isomorphism
$\kf_{\sigma_g}L_g\cong \kf_{\alpha^{-1}}L_g$ described above.

Let $g_1,\ldots,g_t$ be a set of representatives of orbits of $L$ on
$G$.
Write $p_i:=p_{g_i}$, $\sigma_i:=\sigma_{g_i}$, $L_i:=L_{g_i}$,
and $H_i:=H_{g_i}$.
By the above analysis, each simple $H$-module contains
one of $\kf p_1,\ldots,\kf p_t$ on restriction to $(\kf G)^*$.
Thus each arises from a simple $\kf_{\sigma_i}L_i$-module as stated.

\begin{prop}[Kashina-Mason-Montgomery]\label{KMM}
The simple $H$-modules
are precisely the modules $\widehat{V}:=H\otimes_{H_{i}}(p_{i}\otimes V)$ 
induced from $H_{i}$, 
where $V$ is a simple $\kf_{\sigma_{i}}L_{i}$-module, and $i$ ranges over
$\{1,\ldots,t\}$.  
The action
of $H_{i}$ on $V':=p_{i}\otimes V$ is given by 
$$(p_h \overline{x})
(p_{i}\otimes v)=\delta_{g_i,h}p_{i}\otimes\overline{x}\cdot v.$$
\end{prop}

More generally, for any $g\in G$ and $\kf_{\sigma_g}L_g$-module $U$,
we let $\widehat{U}:=H\otimes_{H_g}(p_g\otimes U)$, where $U':=p_g\otimes U$
is the $H_g$-module given by $(p_h \overline{x})(p_g\otimes u)
=\delta_{g,h}p_g\otimes \overline{x}\cdot u$.
As the idempotent $\sum_{y\in L/L_g}p_{ {}^yg}$ acts as the identity
on such a module, and as 0 on $H$-modules corresponding to other orbits of
$L$ on $G$, each $H$-module is a sum of indecomposable modules 
$\widehat{U}$ corresponding to $L$-orbits on $G$.
Each $\widehat{U}$ is then induced from an $H_g$-module of the form
$U'=p_g\otimes U$ where $U$ is a $\kf_{\sigma_g}L_g$-module.
Thus we obtain a decomposition of the category of $H$-modules
into a product of the categories of $\kf_{\sigma_i}L_i$-modules
($i=1,\ldots,t$).
This yields an additive decomposition of Grothendieck groups
\cite[Cor.\ 3.4]{kashina-mason-montgomery02}:
\begin{equation}\label{repring}
K_0(H)\cong \bigoplus_{i=1}^tK_0(\kf_{\sigma_i}L_i).
\end{equation}
As $H$ is a Hopf algebra, the tensor product of modules induces
a ring structure on $K_0(H)$, but the Grothendieck groups
$K_0(\kf_{\sigma_i}L_i)$ are not necessarily rings themselves.
We will see next how the product on $K_0(H)$ behaves with
respect to the additive decomposition (\ref{repring}), 
that is we will give a formula
for the tensor product of two $H$-modules corresponding to 
$\kf_{\sigma_i}L_i$- and $\kf_{\sigma_j}L_j$-modules.

If $H$ is not semisimple, we may alternatively replace $K_0(H)$
by the representation ring of $H$, that is the ring generated by
isomorphism classes of $H$-modules with direct sum for addition and
tensor product for multiplication.
The image of each module in this ring is the sum of images of indecomposable
modules, each corresponding to an $L$-orbit on $G$.
Theorem \ref{tensor} below governs the tensor product of two such
indecomposable modules.

First suppose that $g,h\in G$ and $x,y\in L_g\cap L_h$.
As $\Delta$ is an algebra map, the following relation holds:
\begin{equation}\label{sigmatau}
\sigma_g(x,y)\sigma_h(x,y)=
\sigma_{gh}(x,y)
\tau^{-1}_{g,h}(x)\tau^{-1}_{g,h}(y)
\tau_{g,h}(xy)
\end{equation}
(see \cite[(4.8)]{kashina-mason-montgomery02}).
That is, $\sigma_{gh}$ is cohomologous
to $\sigma_g\cdot \sigma_h$ on $L_g\cap L_h$.
Therefore there is an isomorphism
\begin{equation}\label{psi}
\psi: \kf_{\sigma_{gh}}(L_g\cap L_h) \stackrel{\sim}{\longrightarrow} 
\kf_{\sigma_g\cdot\sigma_h}(L_g\cap L_h)
\end{equation}
given by $\psi(\overline{x})=\tau_{g,h}(x)\overline{x}$.
This will be important in taking tensor products of modules,
as the tensor product of a $\kf_{\sigma_g}(L_g\cap L_h)$-module
with a $\kf_{\sigma_h}(L_g\cap L_h)$-module is naturally a $\kf_{\sigma_g
\cdot\sigma_h}(L_g\cap L_h)$-module.
We will want to induce such a module 
to a $\kf_{\sigma_{gh}}L_{gh}$-module.
This involves first applying the isomorphism (\ref{psi}), and then applying
induction from the subalgebra $\kf_{\sigma_{gh}}(L_g\cap L_h)$
to the algebra $\kf_{\sigma_{gh}}L_{gh}$.

Note that ${}^xL_j=L_{{}^xg_j}$.
We will use an arrow down ($\downarrow$) to denote restriction of modules to
the indicated subalgebra of $\kf_{\sigma_g}L_g$, and an arrow up 
($\uparrow$) to denote tensor induction of
modules from a subalgebra to $\kf_{\sigma_g}L_g$.
The following theorem does not give a formula for decomposing
the tensor product of two simple (respectively, indecomposable)
modules fully into a direct sum of 
{\em simple} (respectively, indecomposable) 
components, although it is a first step towards such a decomposition.

\begin{thm}\label{tensor}
Let $V$ be a $\kf_{\sigma_{i}}L_{i}$-module, $W$ a 
$\kf_{\sigma_{j}}L_{j}$-module, and $\widehat{V}, \widehat{W}$ the
corresponding $H$-modules as described in Proposition \ref{KMM}.
Then as $H$-modules,
$$\widehat{V}\otimes \widehat{W}\cong \sum_{x\in D}
\widehat{U(x)},$$
where $D$ is a set of representatives of double cosets 
$L_i\backslash L/L_j$, and $U(x)$ is the $\kf_{\sigma_{g_i({}^xg_j)}}L_{g_i
({}^xg_j)}$-module
$$U(x)=(V\downarrow^{L_{i}}_{L_{i}\cap {}^xL_{j}}\otimes {}^xW
 \downarrow^{{}^xL_{j}}
 _{L_{i}\cap {}^xL_{j}})\uparrow_{L_{i}\cap {}^xL_{j}}^{L_
{g_i({}^xg_j)}}.$$
\end{thm}

We obtain a formula for
the product in $K_0(H)$ in terms of the additive
decomposition (\ref{repring}) by identifying $\widehat{U(x)}$
with the isomorphic $H$-module
$\widehat{ {}^yU(x)}$ where $k=k(x)$ and $y=y(x)$ are chosen
so that $g_k={}^yg_i {}^{yx}g_j$.
Then the formula in the theorem closely resembles the formulas in
Corollary \ref{product} and Theorem \ref{productformula}.

\begin{proof}[Proof of Theorem \ref{tensor}]
First we will check that the underlying $(\kf G)^*$-modules are isomorphic.
The Grothendieck ring of $(\kf G)^*$-modules is just the group algebra $\Z G$.
The underlying $(\kf G)^*$-module of $\widehat{V}$ (respectively, of 
$\widehat{W}$)
is $\dim V$ (respectively, $\dim W$) copies of the sum ${\mathcal O}(g_i)$
of the elements in the orbit
of $g_i$ (respectively, ${\mathcal O}(g_j)$ of $g_j$).
The underlying $(\kf G)^*$-module of $\widehat{U(x)}$ is 
$|L_{g_i({}^xg_j)}:L_{i}\cap {}^xL_{j}| (\dim V)(\dim W)$ copies of
the sum ${\mathcal O}(g_i({}^xg_j))$ of elements in the orbit
of $g_i({}^xg_j)$.
In $\Z G$, the product of sums of orbit elements is
$${\mathcal O}(g_i)\cdot{\mathcal O}(g_j)=\sum_{x\in D}
|L_{g_i({}^xg_j)}:L_{i}\cap {}^xL_{j}|{\mathcal O}(g_i({}^xg_j)).$$
(This follows from standard properties of the trace map for the
$L$-algebra $\Z G$.
See for example \cite{thevenaz??}.)
Multiplying both sides by $(\dim V)(\dim W)$, we see that
the underlying $(\kf G)^*$-modules in the statement of the theorem
are indeed isomorphic.

We will next identify the appropriate
$H$-submodules of $\widehat{V}\otimes \widehat{W}$,
keeping in mind its underlying $(\kf G)^*$-module structure.
For each $x\in D$, we will describe the $H_{g_i({}^xg_j)}$-module
generated by $p_{g_i}\widehat{V}\otimes p_{{}^xg_j}\widehat{W}$.
Consider the map
\begin{equation}\label{VW}
\left(V\downarrow^{L_i}_{L_i\cap {}^xL_j} \otimes {}^xW\downarrow
  ^{{}^xL_j}_{L_i\cap {}^xL_j}\right)\uparrow^{L_{g_i({}^xg_j)}}
  _{L_i\cap {}^xL_j} \longrightarrow H_{g_i({}^xg_j)}\left(p_{g_i}
\widehat{V}\otimes p_{{}^xg_j}\widehat{W}\right)
\end{equation}
given by $\overline{y}(v\otimes w)\mapsto \overline{y}((p_{g_i}\otimes v)
\otimes
(p_{{}^xg_j}\otimes w))$ for $y$ ranging over a set of coset representatives
of $L_i\cap {}^xL_j$ in $L_{g_i ({}^xg_j)}$.
Note that distinct coset representatives $y$ will generate disjoint
subspaces $\overline{y}(p_{g_i}\widehat{V}\otimes p_{{}^xg_j}\widehat{W})$,
and each has dimension $(\dim V)(\dim W)$, so the above map is bijective.
It may be checked that the action of
$\kf_{\sigma_{g_i({}^xg_j)}}L_{g_i({}^xg_j)}$ on the left side of
(\ref{VW}) corresponds to the action of $H_{g_i({}^xg_j)}$ on the right
side of (\ref{VW}), using the isomorphism (\ref{psi}).
The induced $H$-module then has dimension $|L:L_i\cap {}^xL_j|
(\dim V)(\dim W)$, as desired.
\end{proof}

For each $g\in G$, let $A(g):=K_0(\kf_{\sigma_g}L_g)$, 
the Grothendieck group of $\kf_{\sigma_g}L_g$-modules.
If $H$ is not semisimple, we may alternatively take $A(g)$ to be
the additive group generated by $\kf_{\sigma_g}L_g$-modules, with
direct sum for addition.
Additively, the representation ring of $H$ decomposes into a direct sum
of such groups, analogous to (\ref{repring}).
The following statements apply equally well to this ring.

We obtain an isomorphism $p_g\overline{L_g}
\stackrel{\sim}{\longrightarrow}p_{{}^xg}\overline{L_{{}^xg}}$ of
subalgebras of $H$ for each $x\in L$, via the action of $L$ 
(see (\ref{inner})).
This induces an isomorphism $\kf_{\sigma_g}L_g
\stackrel{\sim}{\longrightarrow} \kf_{\sigma_{{}^xg}}L_{{}^xg}$.
Thus we have conjugation maps $c_{g,x}:A(g)\stackrel{\sim}{\longrightarrow}
A({}^xg)$.
Define $m_{g,h}:A(g)\times A(h)\rightarrow A(gh)$ by the following for a
$\kf _{\sigma_g}L_g$-module $V_g$ and a $\kf _{\sigma_h}L_h$-module $W_h$:
$$ m_{g,h}(V_g,W_h)=(V_g\downarrow^{L_g}_{L_g\cap L_h}\otimes
W_h\downarrow^{L_h}_{L_g\cap L_h})\uparrow_{L_g\cap L_h}^{L_{gh}}.$$
Then properties (H1)--(H3) are straightforward, and (H4$'$)
is equivalent to the product formula of Theorem \ref{tensor} (compare with
Corollary \ref{product}).
As $c_{g,x}=\id$ when $x\in C(g)$, by (\ref{decomposition}) we have
additive isomorphisms
$$A^L\cong \bigoplus_{i=1}^tK_0(\kf_{\sigma_i}L_i)\cong
K_0(H),$$
and in fact the product on $A^L$ given by Theorem \ref{invariants}(ii)
coincides with that on $K_0(H)$.

\end{document}